\documentclass[]{svmult}

\usepackage{makeidx}         
\usepackage{graphicx}        
\usepackage{multicol}        
\usepackage[bottom]{footmisc}
\usepackage{natbib}
\usepackage{url}

\renewcommand{\PackageWarningNoLine}[2]{}
\usepackage{amsmath}
\usepackage{amsfonts}

\usepackage{color}

\begin{document}

\title*{Parareal for diffusion problems with space- and time-dependent coefficients}
\author{Daniel Ruprecht\inst{1} 
\and Robert Speck\inst{1,2} 
\and Rolf Krause\inst{1}
}
\institute{Institute of Computational Science, Universit{\`a} della Svizzera italiana, 
Lugano, Switzerland, \texttt{\{daniel.ruprecht,rolf.krause\}@usi.ch} \and
J\"ulich Supercomputing Centre, Forschungszentrum J\"ulich GmbH, J\"ulich, Germany, \texttt{r.speck@fz-juelich.de}}

%
%
\maketitle


\newcommand{\Tt}[1]{\mathbf{#1}}

\section{Introduction}
The very rapidly increasing number of cores in state-of-the-art supercomputers fuels both need for and interest in novel numerical algorithms inherently designed to feature concurrency.
In addition to the mature field of space-parallel approaches (e.g. domain decomposition techniques), time-parallel methods that allow concurrency along the temporal dimension are now an increasingly active field of research, although first ideas, like in~\cite{Nievergelt1964}, go back several decades.
A prominent and widely studied algorithm in this area is Parareal, introduced in~\cite{LionsEtAl2001}, which has the advantage that one can couple and reuse classical time-stepping schemes in an iterative fashion to parallelize in time.
However, there also exist a number of other approaches, e.g. the "parallel implicit time algorithm" (PITA) from~\cite{FarhatEtAl2003}, the "parallel full approximation scheme in space and time" (PFASST) from~\cite{EmmettMinion2012} or "revisionist integral deferred corrections" (RIDC) from~\cite{ChristliebEtAl2010} to name a few.
Parareal in particular and temporal parallelism in general has been considered early as an addition to spatial parallelism in order to extend strong scaling limits, see~\cite{MadayTurinici2005}.
Efficacy of this approach in large-scale parallel simulations on hundred of thousands of cores has been demonstrated for the PFASST algorithm in~\cite{SpeckEtAl2012}.

For Parareal, multiple works exist that demonstrate its efficiency for diffusion problems:
\cite{GanderVandewalle2007_SISC} prove super-linear convergence of Parareal for the standard 1D heat equation.
A  more general theorem showing super-linear convergence for nonlinear ODEs is proven by~\cite{GanderHairer2008}, while~\cite{Bal2005} presents a convergence theorem for linear parabolic PDEs with constant coefficients.
The present paper investigates the effect of space- and time-dependent coefficients in the two-dimensional heat equation on the convergence of Parareal.
This is done by means of numerical examples, including one that shows how convergence of Parareal can be estimated by the maximum singular value of a Parareal iteration matrix.
\section{Parareal}
To match the numerical examples in \S\ref{sec:examples}, the presentation of Parareal given here starts with an initial value problem
\begin{equation}
	M y_{t}(t) = f(y(t), t), \ y(0) = b \in \mathbb{R}^{d}, \ t \in [0,T],
\end{equation}
with a mass matrix $M$ and right-hand side $f$ arising from a finite element discretization of a partial differential equation.
Let $(t_n)_{n=0}^{N}$ with $t_0 = 0$ and $t_{N} = T$ be a decomposition of $[0,T]$ into $N$ so-called time-slices $[t_{n}, t_{n+1}]$ which, for the sake of simplicity, are assumed to be of equal length here.
Furthermore, let $y_{n}$ be an approximation to the solution at $t_{n}$, that is \linebreak $y_{n} \approx y(t_{n})$.\par
Denote by $\mathcal{F}$ a "fine", computationally expensive and accurate integration method with a time step $\delta t$ (e.g. a higher-order Runge-Kutta method) and by $\mathcal{G}$ a "coarse", computationally cheap and probably inaccurate method with a time step $\Delta t \gg \delta t$ (e.g. implicit Euler).
Assume here that the constant length of the time-slices is a multiple of both $\delta t$ and $\Delta t$, so that fine as well as coarse method can integrate over one time-slice using a fixed integer number of time-steps.
Denote the result of integrating over the slice $[t_{n}, t_{n+1}]$, starting from an initial value $y$ at $t_{n}$, using the fine or coarse method as $\mathcal{F}(y, t_{n+1}, t_{n})$ and $\mathcal{G}(y, t_{n+1}, t_{n})$ respectively.
Serial integration using the fine method would then correspond to computing
\begin{equation}
	\label{eq:fine}
	y_{n+1} = \mathcal{F}(y_{n}, t_{n+1}, t_{n}), \ n=0,\ldots,N-1,
\end{equation}
step-by-step with $y_0 := b$.
Instead, Parareal computes the iteration given by
\begin{equation}	
	\label{eq:parareal}
	y^{k+1}_{n+1} = \mathcal{G}(y^{k+1}_{n}, t_{n+1}, t_{n}) + \mathcal{F}(y^{k}_{n}, t_{n+1}, t_{n}) - \mathcal{G}(y^{k}_{n}, t_{n+1}, t_{n})
\end{equation}
where the evaluation of the fine method over the $N$ time-slices can be distributed over $N$ processors (see~\cite{LionsEtAl2001} for details).
The iteration converges to the serial fine solution as $k \to N$. Speedup can be achieved if $\mathcal{G}$ is sufficiently cheap compared to $\mathcal{F}$ and if the iteration converges in $K \ll N$ iterations.
Therefore, rapid convergence is critical for Parareal to be efficient.
In the examples below, the defect 
\begin{equation}
	d^{k} := \max_{i=0,\ldots,N} \left\| y_{i} - y^{k}_{i} \right\|_{\infty}
\end{equation}
between the solution provided by the Parareal iteration~\eqref{eq:parareal} after $k$ iterations and the serial fine solution~\eqref{eq:fine} is used to measure convergence.
\section{Heat equation with non-constant coefficients}\label{sec:examples}
The test problem used here to study the convergence of Parareal for non-constant coefficients is the two-dimensional heat equation
\begin{equation}
	u_{t}(x,y,t) = \nu(t) \nabla \cdot \left( a(x,y) \nabla u(x,y,t) \right)
\end{equation}
 on a square $\Omega = [0,1]^{2}$. 
 The initial values are given by
 \begin{equation}
	u_0(x,y) =  \exp\left[ - \left( (x-0.5)^2 + (y-0.5)^2 \right) / \sigma^{2} \right], \ \sigma=0.35,
\end{equation}
and the problem is run until $T = 4.0$.
The interval $[0,T]$ is divided into $N = 40$ time-slices and an implicit Euler method with $\Delta t = 1/100$ is used for $\mathcal{G}$ and a third order RadauIIA(3) method with $\delta t = 1/200$ for $\mathcal{F}$.
The spatial domain $\Omega$ is divided into three "strips" 
\begin{align}
	\Omega_{1} &= [0, x_0) \times [0, 1], \\
	\Omega_{2} &= [x_0, x_0 + w) \times [0,1], \\
	\Omega_{3} &= [x_0 + w, 1] \times [0,1],
\end{align}
and a different constant value for $a$ is prescribed on every strip, that is
\begin{equation}
	a(x,y) = \left\{ \begin{array}{c @{ \ : \ }l} a_1 & (x,y) \in \Omega_{1} \\ a_2 & (x,y) \in \Omega_2 \\ a_3 & (x,y) \in \Omega_3. \end{array} \right.
\end{equation}
Furthermore, the effect of varying the width $w$ of the middle strip $\Omega_2$ is investigated.
Conforming triangle meshes aligned with the strips $\Omega_i$ are generated for values of $w \in \left\{ 0.2, 0.1, 0.05, 0.02 \right\}$.
Then, for every value of $w$, a number of uniform refinement steps is performed in order to produce meshes of comparable mesh width.
After refinement, the minimum element sizes for the different values of $w$ range from $h_{\rm min} = 0.01$ to $h_{\rm min} = 0.005$ and the maximum element sizes from $h_{\rm max} = 0.02$ to $h_{\rm max} = 0.035$, so that the resolutions are comparable.
All experiments reported below use linear finite elements, but preliminary tests not documented here suggest that the results are not significantly affected by the use of higher-order FEM.
Homogeneous Dirichlet boundary conditions are employed.
Simulations are run with $a_1 = a_3 = 0.01$ fixed and $a_2 \in \left\{ 0.01, 1.0, 100 \right\}$, resulting in ratios $\Delta a = a_2 / a_3 = a_2 / a_1 \in \left\{ 1, 100, 10000 \right\}$.
\subsection{Space-dependent coefficients}
\begin{figure}[th]
	\centering
	\begin{minipage}{0.45\textwidth}
	\centering
	\hspace{3em} Strip width $w = 0.2$
	\includegraphics[width=0.99\textwidth]{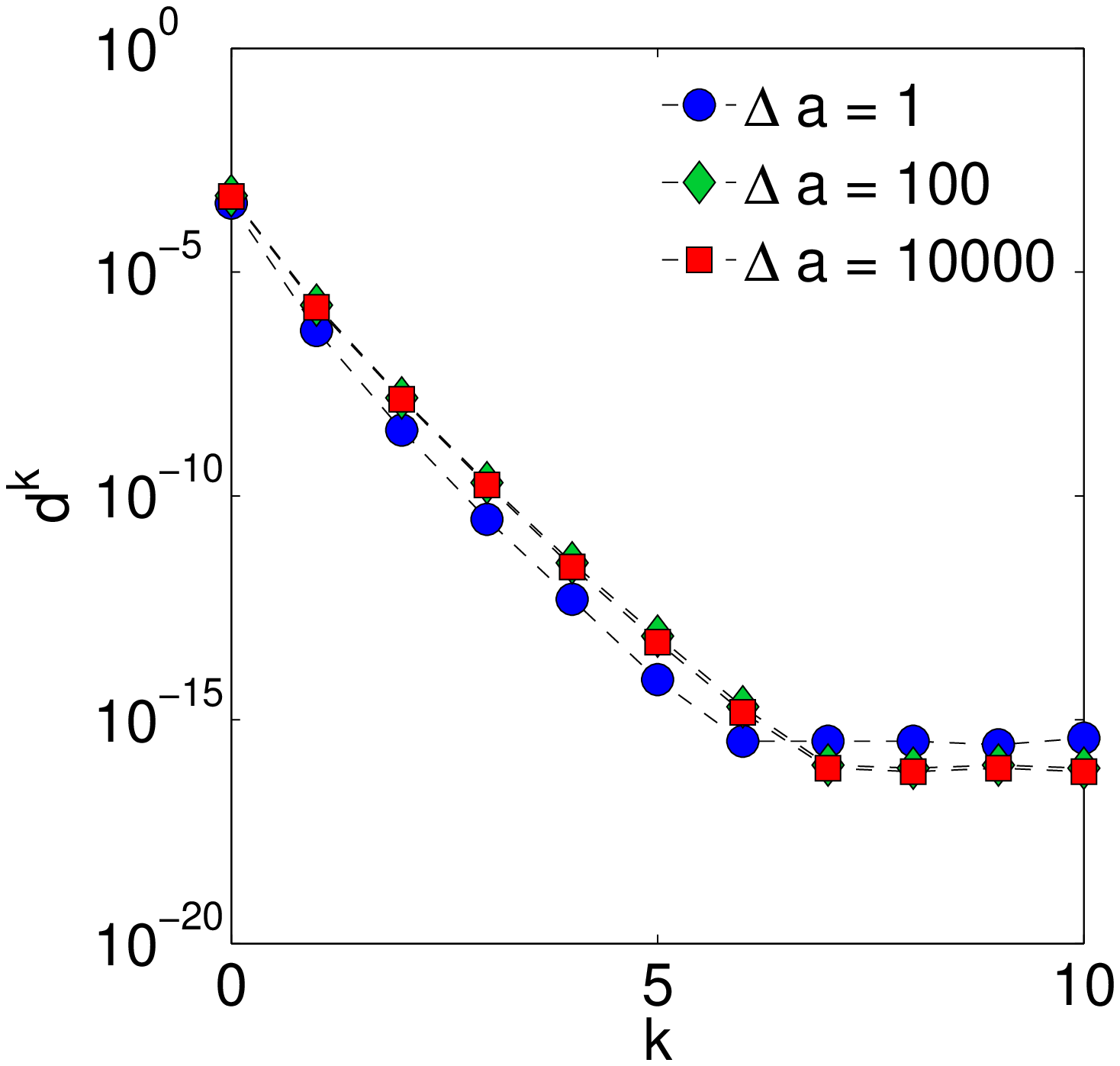}
	\end{minipage}
	\begin{minipage}{0.45\textwidth}
	\centering
	\hspace{3em} Strip width $w = 0.02$
	\includegraphics[width=0.99\textwidth]{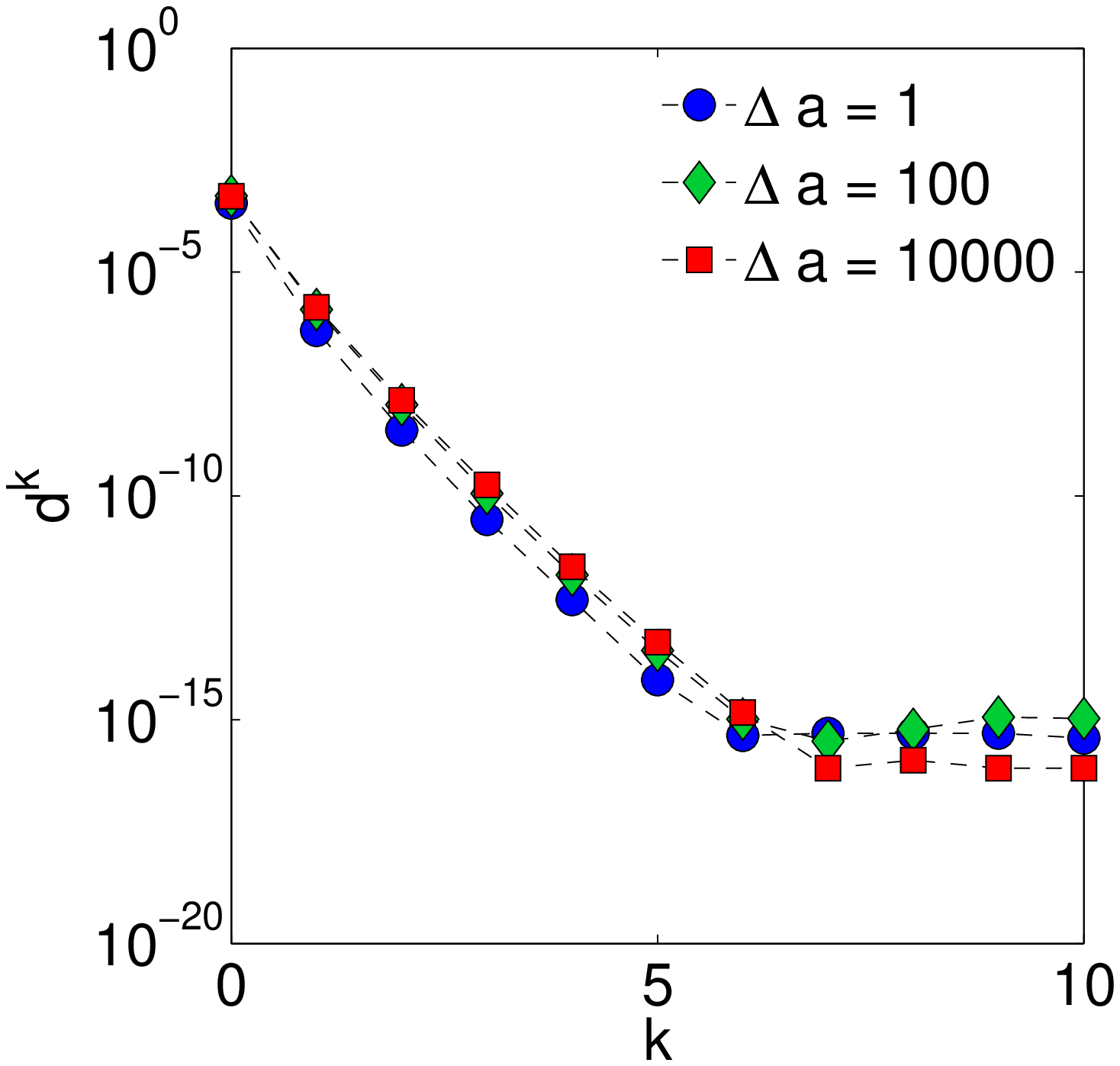}
	\end{minipage}	
	\caption{Defect $d^{k}$ between Parareal and the serial fine solution versus the iteration number $k$ depending on the magnitude of the jump in the diffusion coefficient from $\Omega_{1}$, $\Omega_{3}$ to $\Omega_{2}$.}\label{fig:space_coeff}
\end{figure}

First, set $\nu \equiv 1$ in order to study only the effect of spatially varying coefficients.
Figure~\ref{fig:space_coeff} shows the resulting convergence of Parareal for the different values of $\Delta a$ and $w=0.2$ (left) and $w=0.02$ (right).
Convergence in the cases with jumping coefficients is slightly slower, but the effect is very small. 
Also, the reduction in convergence speed seems to be rather independent of the magnitude of the jump in the diffusion coefficient: In both plots, the lines for $\Delta a = 100$ and $\Delta a = 10000$ are more or less indistinguishable.\par
Convergence of Parareal is utterly oblivious against the width $w$ of the middle strip $\Omega_2$: When plotting the defects for fixed $\Delta a$ and different values of $w$, the resulting data points all essentially coincide so that the corresponding plots are rather uninteresting and are therefore omitted.
\subsection{Space- and time-dependent coefficients}
\begin{figure}[th]
	\centering
	\begin{minipage}{0.45\textwidth}
		\centering
		\hspace{3em} Slow change of $\nu$: $\alpha=1$
		\includegraphics[width=0.99\textwidth]{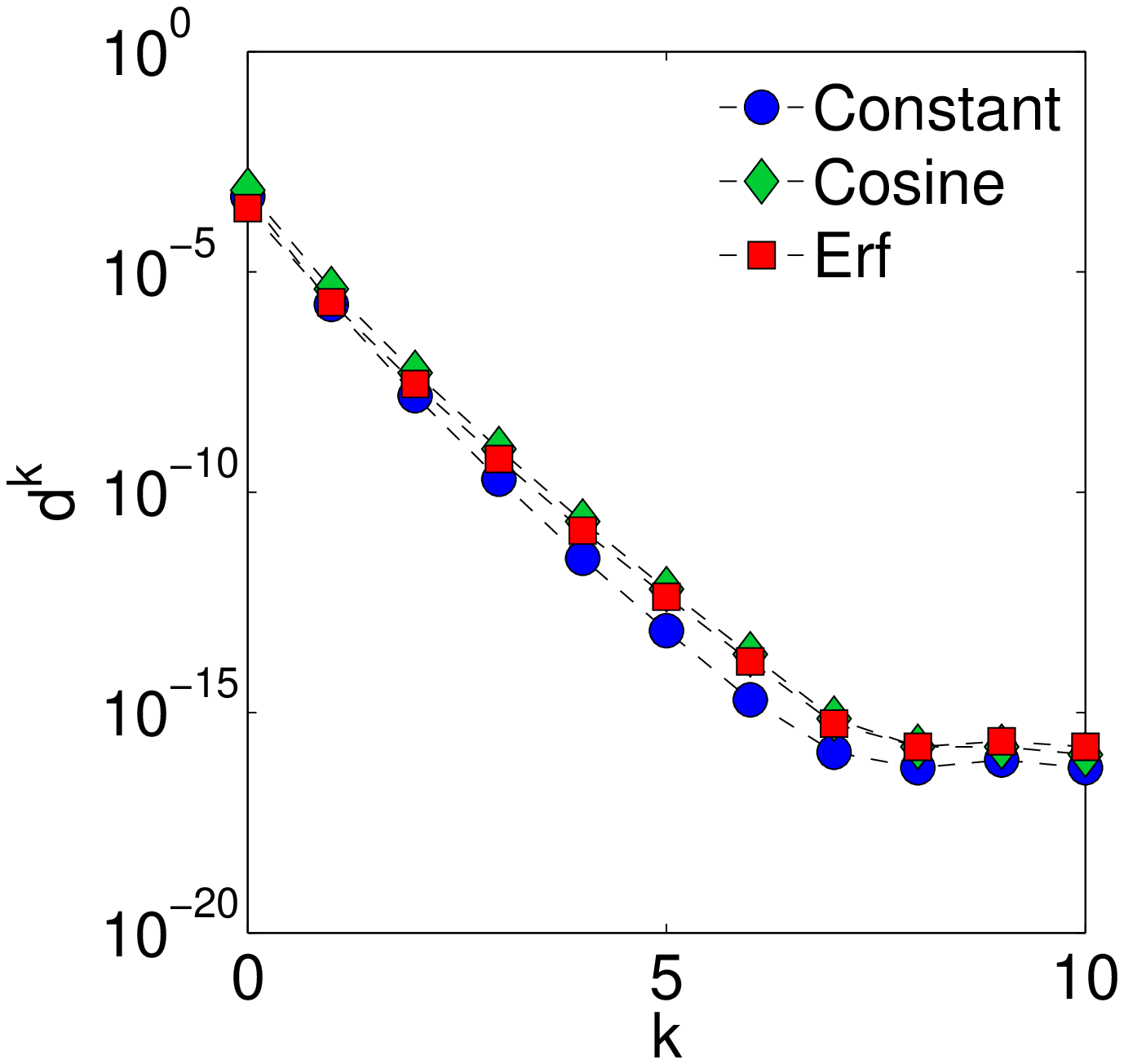}
	\end{minipage}
	\begin{minipage}{0.45\textwidth}
		\centering
		\hspace{3em} Fast change of $\nu$: $\alpha=10$
		\includegraphics[width=0.99\textwidth]{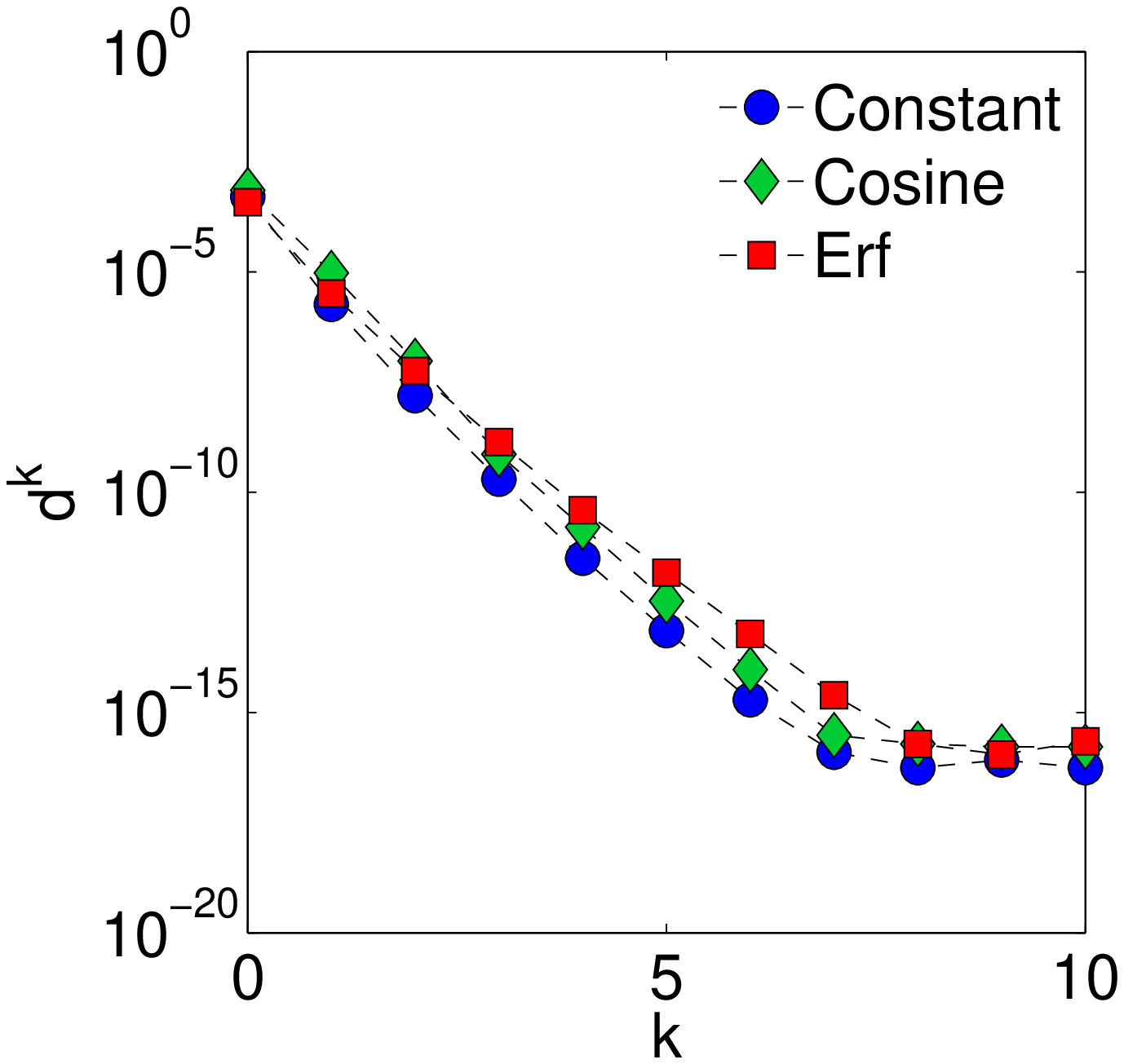}
	\end{minipage}	
	\caption{Defect $d^{k}$ of Parareal versus the iteration number $k$ for different time-dependent $\nu$-profiles with $\Delta a = 100$, $\alpha = 1$ (left) and $\alpha = 10$ (right).}\label{fig:time_coeff}
\end{figure}

To investigate the effect of a time-dependent diffusion coefficient on the convergence of Parareal, fix the strip width to $w=0.2$ and the coefficient jump to $\Delta a = 100$.
Furthermore, use the following three different profiles for the time-dependent diffusion coefficient $\nu$ :
\begin{align}
	\nu(t) &= 1 \quad &\text{("constant")}, \\
	\nu(t) &= \frac{1}{2}\left(1 + \cos\left( \alpha \frac{\pi}{2} t \right) \right) \quad &\text{("cosine")}, \\
	\nu(t) &= \frac{1}{2} \left(1 + \textrm{erf}( \alpha (t-2) ) \right) \quad &\text{("erf")}. \label{eq:err_f_nu}
\end{align}
Initial value and boundary conditions are set as described above. 
Two sets of simulations are performed, one with $\alpha = 1$ corresponding to a very slowly changing $\nu$ and one with $\alpha=10$ corresponding to a more rapid change.
The resulting convergence of Parareal is shown in Figure~\ref{fig:time_coeff}.
In both cases, the slow as well as the fast varying one, Parareal's convergence is only marginally affected by the space- and time-dependent diffusion coefficients.
The resulting defects are slightly larger than for the reference case and the difference is a little more pronounced for $\alpha = 10$, but the overall effect is not drastic:
In the fast varying case with the error function profile~\eqref{eq:err_f_nu}, Parareal requires only a single additional iteration compared to the constant reference in order to reach the same defect level.
\subsection{Error bound from singular values}
\begin{figure}[th]
	\centering
	\begin{minipage}{0.45\textwidth}
		\centering
		\hspace{3em} $\Delta a = 1$
		\includegraphics[width=0.99\textwidth]{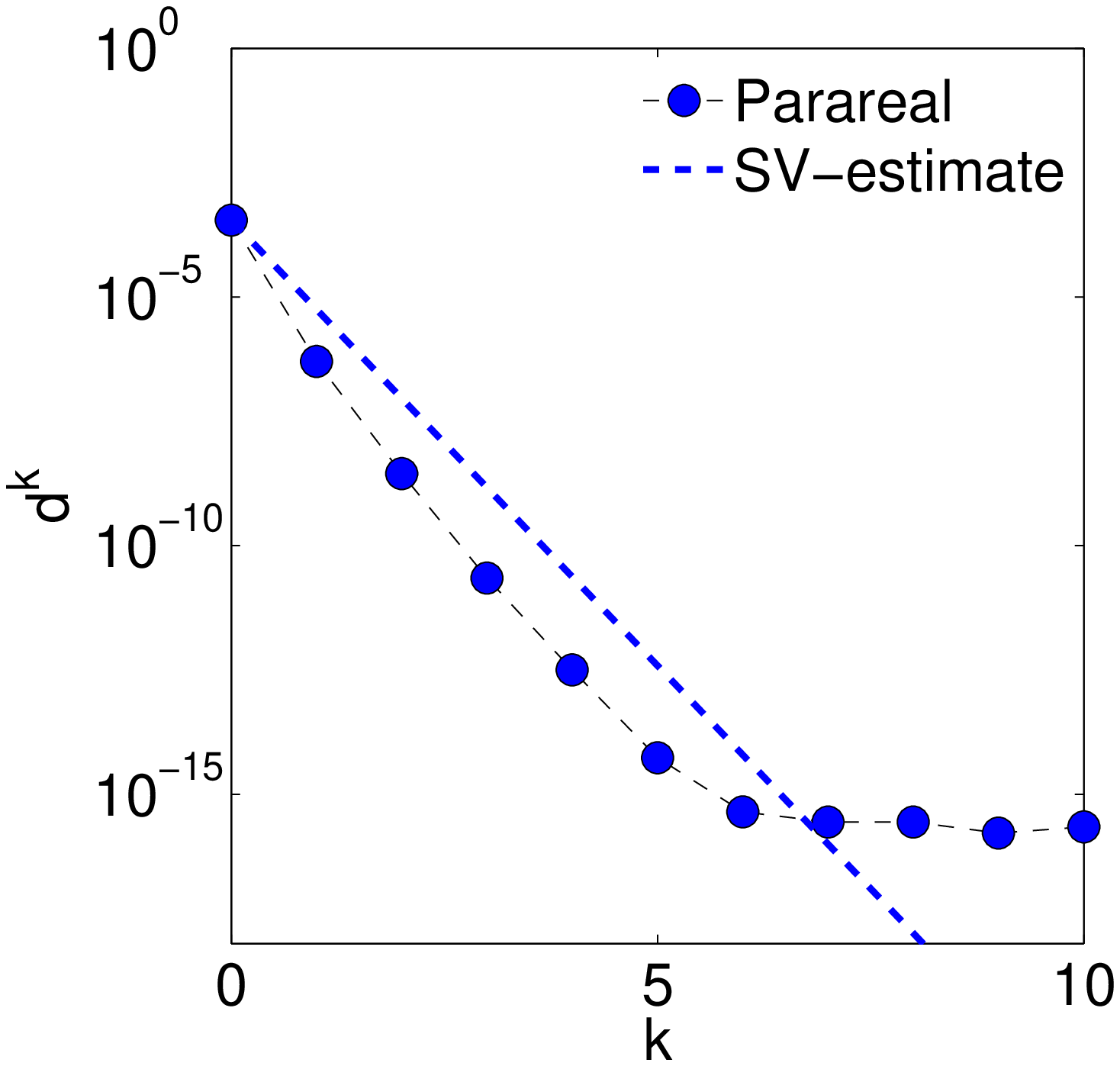}
	\end{minipage}
	\begin{minipage}{0.45\textwidth}
		\centering
		\hspace{3em} $\Delta a = 10000$
		\includegraphics[width=0.99\textwidth]{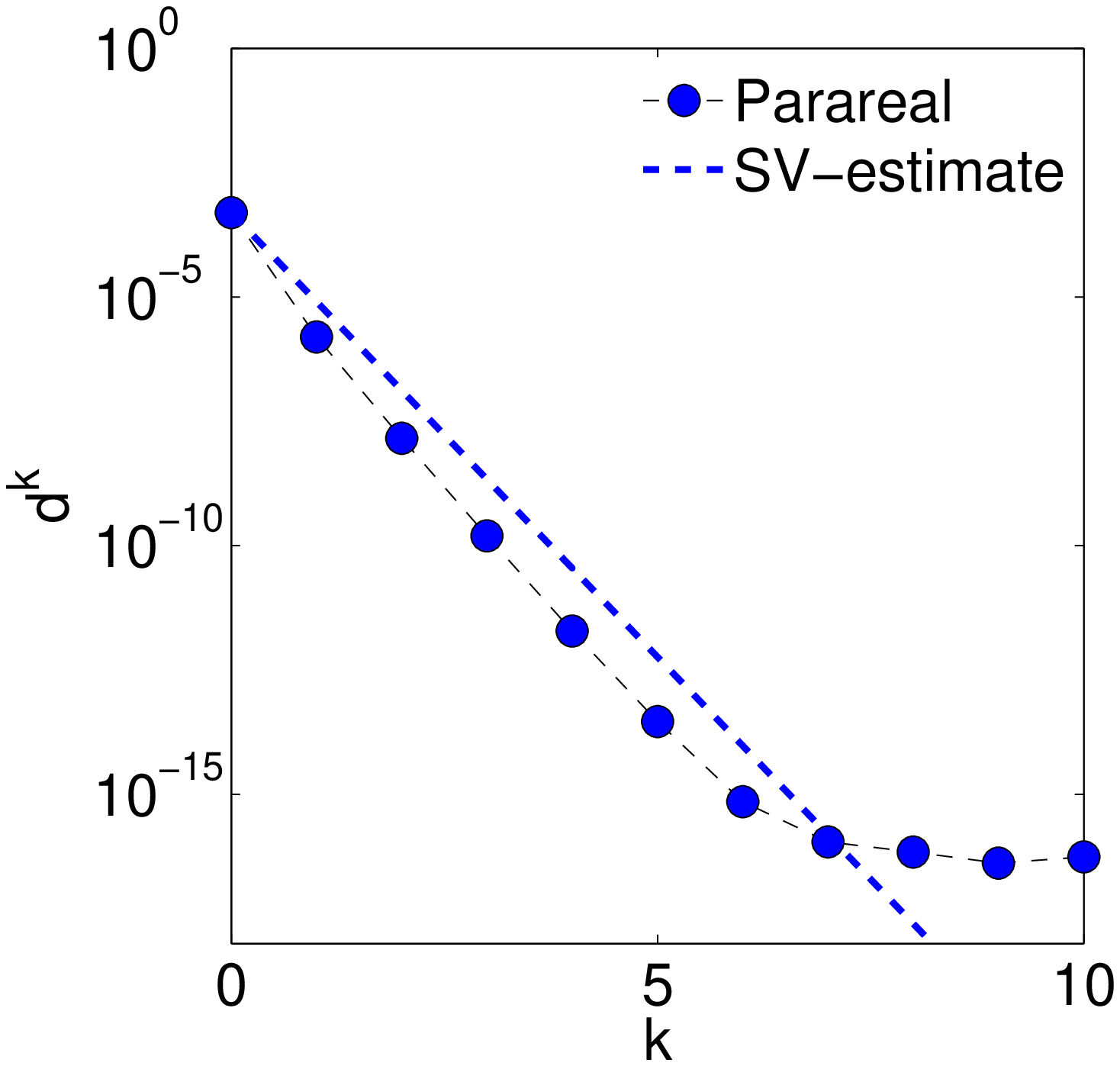}		
	\end{minipage}	
	\caption{Convergence of Parareal and error estimate from the largest singular value of the Parareal iteration matrix (dashed line).}\label{fig:conv_svd}
\end{figure}

Parareal can be considered as a fixed point iteration, see e.g.~\cite{AmodioBrugnano2009} or~\cite{FriedhoffEtAl2013} for more detailed explanations.
For $\nu \equiv 1$ and the linear problem considered here, the action of the propagators $\mathcal{F}$ and $\mathcal{G}$ can be expressed as multiplication by matrices $G$ or $F$. 
Then, running the fine or coarse method over all $N$ time-slices can be expressed as inversion of size $Nd \times Nd$ matrices
\begin{equation}
	\Tt{M}_{f} = \begin{bmatrix} I & \ldots \\ -F & I &  \\  & \ddots & \ddots & \\ & & & -F & I \end{bmatrix}, \quad \Tt{M}_{g} = \begin{bmatrix} I & \ldots \\ -G & I &  \\  & \ddots & \ddots & \\ & & & -G & I \end{bmatrix},
\end{equation}
so that computing the fine solution through~\eqref{eq:fine} corresponds to a block-wise solution of $\Tt{M}_{\rm f} \Tt{y} = \Tt{b}$ with $\Tt{y}= ( y_0, \ldots, y_{N})^{\rm T}$ and $\Tt{b} = ( b, 0, \ldots, 0)^{\rm T}$.
The Parareal iteration~\eqref{eq:parareal} can then be written as the preconditioned fixed point iteration
\begin{equation}
	\Tt{M}_{g} \Tt{y}^{k+1} = \left( \Tt{M}_{g} - \Tt{M}_{f} \right) \Tt{y}^{k} + \Tt{b},
\end{equation}
where inverting $\Tt{M}_{g}$ corresponds to running the coarse method.
A straightforward computation shows that the iteration matrix $\Tt{I} - \Tt{M}_{g}^{-1} \Tt{M}_{f}$ is nilpotent and thus that its spectral radius is zero, corresponding to the well-known fact that Parareal always converges to the fine solution after $N$ iterations, see e.g.~\cite{GanderVandewalle2007_SISC} (although Parareal won't provide any speedup in this case).
A bound for the convergence rate can be obtained by computing the maximum singular value instead.
In order to keep the size of the iteration matrix manageable, the example studied here uses only the $w=0.2$ geometry with a coarser grid with $h_{\rm min} = 0.04$, $h_{\rm max} = 0.068$ and only $N = 20$ time-slices.
The maximum singular values of the iteration matrix are computed with Matlab's {\sc svds} function and are $\sigma_{\rm max} \approx 0.162$ for $\Delta a = 1$ and $\sigma_{\rm max} \approx 0.163$ for $\Delta a = 10000$: The minimal difference gives an additional indication that the coefficient jump should not influence Parareal's convergence.
Figure~\ref{fig:conv_svd} shows the convergence rates of Parareal for this example for $\Delta a = 1$, i.e.~with a constant coefficient (left) and with $\Delta a  = 10000$ (right), as well as the estimate $d^{0} \times (\sigma_{\rm max})^{k}$ resulting from the maximum singular value.
In both cases, actual convergence is a little better than expected but $\sigma_{\rm max}$ gives a reasonable estimate.
Again, the jumping coefficients affect Parareal's convergence only marginally.
Note that interpreting variants like the "Krylov-subspace-enhanced Parareal", introduced in~\cite{GanderPetcu2008} and studied further in~\cite{RuprechtKrause2012}, as a non-stationary fixed point iteration could be an interesting approach for a mathematical analysis.
\section{Conclusions}
The paper presents a numerical study of the convergence behaviour of the time-parallel Parareal method for the heat equation with space- and time-dependent coefficients.
It demonstrates that the good convergence of Parareal for diffusive problems is only marginally affected by both jumps in the diffusion coefficients and a diffusion coefficient that changes in time.
For linear problems, Parareal can be interpreted as a preconditioned fixed point iteration and, at least for small enough problems, the iteration matrix and its maximum singular value can be computed numerically.
An example is shown that demonstrates that the largest singular value gives a reasonable estimate for the convergence of Parareal.
Confirming the robustness of Parareal's convergence for three-dimensional diffusion problems with a complicated geometry would be one interesting direction of future research.
\subsubsection*{Acknowledgments.}
This work was supported by the Swiss National Science Foundation (SNSF) under the lead agency agreement through the project "ExaSolvers" within the Priority Programme 1648 ``Software for Exascale Computing'' (SPPEXA) of the Deutsche Forschungsgemeinschaft (DFG). The authors thankfully acknowledge support from Achim Sch\"adle, who provided parts of the used code.
%
%

%
%
\bibliographystyle{plainnat}
\bibliography{../../bibtex/Pint,../../bibtex/Pint_Self}



\end{document}